\documentclass{article}

\usepackage{amsthm,amssymb,color}
\usepackage{geometry}
\usepackage{comment}
\usepackage{amsmath}
\usepackage{graphicx}
\usepackage{amsfonts}
\usepackage{cite}
\usepackage{hyperref}
\usepackage{mathrsfs}
\usepackage{multirow}
\usepackage{listings}
\usepackage{authblk}
\numberwithin{figure}{section}
\numberwithin{table}{section}

\begin{document}
	\title{On the Uniqueness of Functions that Maximize the Crouzeix Ratio}
	\date{\today}
	\author{Kenan Li\footnote{University of Washington, Applied Mathematics Dept., Box 353925, Seattle, WA 98195.}}
	\maketitle
	\begin{abstract}
		Let $A$ be an $n$ by $n$ matrix with numerical range $W(A) := \{ q^{*}Aq : q \in \mathbb{C}^n , ~\| q \|_2 = 1 \}$.  We are interested in functions $\hat{f}$ that maximize $\| f(A) \|_2$ (the matrix norm induced by the vector 2-norm) over all functions $f$ that are analytic in the interior of $W(A)$ and continuous on the boundary and satisfy $\max_{z \in W(A)} | f(z) | \leq 1$. It is known that there are functions $\hat{f}$ that achieve this maximum and that such functions are of the form $B\circ\phi$, where $\phi$ is any conformal mapping from the interior of $W(A)$ to the unit disk $\mathbb{D}$, extended to be continuous on the boundary of $W(A)$, and $B$ is a Blaschke product of degree at most $n-1$. It is not known if a function $\hat{f}$ that achieves this maximum is unique, up to multiplication by a scalar of modulus one. We show that this is the case when $A$ is a $2\times 2$ nonnormal matrix or a Jordan block, but we give examples of some $3\times 3$ matrices with elliptic numerical range for which two different functions $\hat{f}$, involving the same conformal mapping but Blaschke products of different degrees, achieve the same maximal value of $||f(A)||_2$.
	\end{abstract}
	
	\section{Introduction} 
	Let $A$ be an $n$ by $n$ matrix with \textit{numerical range} $W(A) := \{ q^{*}Aq : q \in \mathbb{C}^n , ~\| q \|_2 = 1 \}$. M. Crouzeix conjectured that the ratio,
	\begin{equation}\label{ratio}
		\sup_{f}\dfrac{||f(A)||_2}{\max_{z\in W(A)}|f(z)|},
	\end{equation}
	referred to here as the \textit{Crouzeix ratio}, is always less than or equal to 2\cite{CROUZEIX2007668}, and he proved that it is always less than or equal to 11.08. Later M. Crouzeix and C. Palencia  reduced this bound to $1+\sqrt{2}$\cite{doi:10.1137/17M1116672}, but the conjectured bound of 2 still has not been established. The supremum in (\ref{ratio}) is over all functions $f$ analytic in the interior of $W(A)$ and continuous on the boundary, and $||\cdot||_2$ denotes the matrix norm induced by the vector 2-norm; that is, the largest singular value. Many numerical experiments support Crouzeix's conjecture. See \cite{GREENBAUM2018225}, for instance, where a numerical optimization code was used to try to minimize the inverse of the ratio in (\ref{ratio}).
	
	It is known from Pick-Nevanlinna theory\cite{Crouzeix2004} that there are \textit{extremal} functions $\hat{f}$ that attain the supremum in (\ref{ratio}) and that such functions are of the form $ \mu B\circ\phi$, where $\mu$ is a nonzero scalar, $\phi$ is any conformal mapping from the interior of $W(A)$ to the unit disk $\mathbb{D}$, extended to be continuous on the boundary, and $B$ is a Blaschke product of degree at most $n-1$:
	\begin{equation}
		B(z)=\exp(i\gamma)\prod_{j=1}^{m}\dfrac{z-\alpha_j}{1-\bar{\alpha_j}z}, ~~m \leq n-1,~ |\alpha_j| < 1.
	\end{equation}
	Without loss of generality, we may assume that $\hat{f}$ has ${\cal H}^{\infty}$-norm $1$ on $W(A)$ and then $\mu = 1$ and $\hat{f}$ is a function of the form $B \circ \phi$ that satisfies
	\begin{equation}
		\hat{f} = \arg \max \{ \| f(A) \|_2 : f \in {\cal H}^{\infty} (W(A)) ,~\| f \|_{\infty} = 1 \} .
		\label{fhat}
	\end{equation}
	It is not known if a function $\hat{f}$ that satisfies (\ref{fhat}) is unique, up to multiplication by a scalar of modulus 1. Note that if $\hat{f}:=\hat{B}\circ\phi$ is one such function and if $\psi$ is a different conformal mapping from the interior of $W(A)$ to $\mathbb{D}$ (mapping a different point to the origin), then $\phi=B_1\circ\psi$ for some degree 1 Blaschke product $B_1$, so that $\hat{f}$ can be written in the form ($\hat{B}\circ B_1$)$\circ\psi$, where $\hat{B}\circ B_1$ is another Blaschke product of the same degree as $\hat{B}$. Thus, to study the uniqueness of $\hat{f}$, we can fix a conformal mapping $\phi$ and ask whether the Blaschke product $\hat{B}$ associated with $\phi$ is unique (where we always mean up to multiplication by a scalar of modulus $1$).
	
	In the following sections we show that  for $2\times 2$ matrices and for Jordan blocks, an \textit{extremal Blaschke product} (i.e., a Blaschke product $\hat{B}$ such that $\hat{f}:=\hat{B}\circ\phi$ satisfies (\ref{fhat}), for a given conformal mapping $\phi$) is unique, but we give examples of some $3\times 3$ matrices with elliptical numerical range for which two different Blaschke products of different degrees both satisfy (\ref{fhat}). Since we show a range of parameters for which the extremal Blaschke product has degree $1$ instead of $2$, this result supports numerical observations that the extremal Blaschke product $\hat{B}$ sometimes appears to have degree strictly less than $n-1$ \cite{Discussions}.
	
	\section{Examples of Uniqueness}
	\subsection{2 by 2 Nonnormal Matrices}
	For any 2 by 2 matrix $W(A)$ is an ellipse, and Crouzeix \cite{Crouzeix2004} proved that an exact expression for the supremum in (\ref{ratio}) is
	$$C(W(A),2)=2\exp\left(-\displaystyle\sum_{n\ge 1}\dfrac{(-1)^{n+1}}{n}\dfrac{2}{1+\rho^{n}}\right),$$ where
	$$~\rho=\dfrac{1+\sqrt{1-\epsilon^2}}{\epsilon},$$
	and $\epsilon$ is eccentricity of the ellipse $W(A)$. Here we give an equivalent expression for the bound $C(W(A),2)$ and prove that $\hat{B}$ (or $\hat{f}$) is unique. 
	
\begin{figure}
		\centering
		\includegraphics[width=12cm]{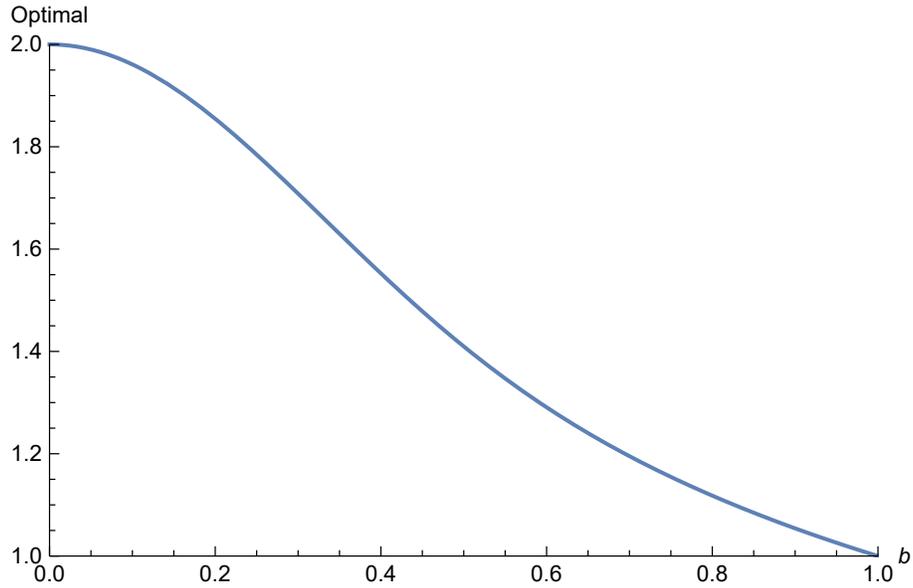}
		\caption{$||\hat{B}(\phi(A))||_2=\sqrt{\frac{k}{b}}$ vs. $b$}
		\label{fig1.3.1}
	\end{figure}

	It can be shown that, after a proper translation and rotation, any $2\times 2$ matrix $A\ne cI$ (trivial case) is unitarily similar to a matrix of the form
	$$\left[\begin{matrix}0 &a\\ d &0 \end{matrix}\right],$$
	where $a>0,d \ge 0$. If $a<d$, it is equivalent to consider $A^T$, so without loss of generality we can assume $a\ge d\ge 0$. The matrix $A$ is normal if and only if $a=d$, and in this case $||f(A)||_2\le\max_{z\in W(A)}|f(z)|$; otherwise, $W(A)$ is an ellipse. G. Szeg\"{o}\cite{Szego} showed that the conformal mapping of the interior of an ellipse with foci points $\pm 1$ onto $\mathbb{D}$, with the origin mapping to itself, is
	\begin{equation}\label{confm}
		\phi_0(z)=\sqrt{k}~\text{sn}\left(\dfrac{2K}{\pi}\sin^{-1}z\right) .
	\end{equation}
	Here $\text{sn}(z)$ is one kind of Jacobi elliptical function, $k$ is the elliptical modulus and $K$ is the quarter period.
	We can multiply $A$ by the scalar $\sqrt{k/(ad)}$ and write
	$$A=\left[\begin{matrix}
	0 &\sqrt{\dfrac{k}{b}}\\
	\sqrt{kb}  &0\\
	\end{matrix}\right],$$
	where $0\le b=\frac{d}{a}\le 1$ and $0 < k= k( b^2 ) < 1$, and then, the conformal mapping maps $A$ to itself, since it maps the eigenvalues of $A$, $\pm\sqrt{k}$, to themselves.
	$$\phi(A)=A=\left[\begin{matrix}
	1 &0 \\
	0 &\sqrt{\dfrac{k}{b}}
	\end{matrix}\right]^{-1}\left[\begin{matrix}
	0 &1\\
	k &0
	\end{matrix}\right]\left[\begin{matrix}
	1 &0\\
	0 &\sqrt{kb} 
	\end{matrix}\right]=V^{-1}TV.$$
	Since $||T||_2=1$ and $\kappa(V)=\sqrt{k/b}$, for any Blaschke product $B$
	$$||B(\phi(A))||_2\le \sqrt{\dfrac{k}{b}}$$
	by von Neumann's inequality\cite{Neumann}, and the equality holds when
	$$\hat{B}(z)=z.$$
	See Figure \ref{fig1.3.1}. 
	
	Further, the extremal $B$ is unique if $A$ is nonnormal ($a>d$ or $b<1$). Suppose an extremal $B$ is of the form
	$$B(z)=\dfrac{z - \alpha} {1 - \bar{\alpha} z},~|\alpha|\le 1,$$ 
	We will show that $\alpha$ must be 0.
	
	From Theorem 5.1 in \cite{caldwell2017extensions}, an extremal Blaschke product $B$ satisfies
	$$B(\phi(A))v_1=||B(\phi(A))||_2u_1,$$
	and
	$$u_1^*v_1=0,$$
	where $u_1$ and $v_1$ are the first left and right singular vectors of  $B(\phi(A))$. Taking
	$$V=[v_1 ~v_2]=[v_1 ~u_1],$$
	we have
	$$V^*V=VV^*=I,$$
	and the diagonal elements of $V^*B(\phi(A))V$ are $0$, i.e., the eigenvalues of $V^*B(\phi(A))V$ or $B(\phi(A))$ are of the form $\pm\lambda$; i.e.,
	$$\dfrac{\sqrt{k}-\alpha}{1-\bar{\alpha}\sqrt{k}}=-\dfrac{-\sqrt{k}-\alpha}{1+\bar{\alpha}\sqrt{k}},~0<{k}\le 1,$$
	which gives us
	$$\alpha = \bar{\alpha}k.$$
	When $0<b<1$, we know $0<k<1$ and $$\alpha=0;$$
	when $b=1$, we know $k=1$ and matrix A is normal.
	
	\subsection{Jordan Blocks}
	We know that for the $n$ by $n$ Jordan block
	$$J = \left[\begin{matrix}
	0 &1 &0 &\cdots &0 \\ 0 &0 &1 &\cdots &0 \\ 
	0 &0 &0 &\ddots &0 \\ 0 &0 &0 &\cdots  &1\\
	0 &0 &0 &\cdots &0
	\end{matrix}\right]_{n\times n},$$
	$W(J)$ is a disk centered at the origin of radius $\cos(\frac{\pi}{n+1})$, and the conformal mapping $\phi$ that maps $W(J)$ to the unit disk $D$ is
	$$\phi(z)=\dfrac{1}{\cos(\frac{\pi}{n+1})}z=cz,~c=\dfrac{1}{\cos(\frac{\pi}{n+1})}>1.$$
	We can rewrite 
	$$\phi(J)=D^{-1}JD,$$
	where
	$$D=\left[\begin{matrix}
	1 & & & \\  &c & & \\ 
	& &\ddots & \\  & & &c^{n-1} 
	\end{matrix}\right]_{n\times n}.$$
	Since $||J||_2=1$, by von Neumann's inequality, 
	$$||B(J)||_2\le 1,$$
	and we know for any Blaschke product $B$
	$$B(\phi(J))=D^{-1}B(J)D,$$
	and 
	$$||B(\phi(J))||_2\le \kappa(D)||B(J)||_2=c^{n-1}||B(J)||_2\le c^{n-1}.$$
	
	One extremal Blaschke product is 
	$$\hat{B}(z)=z^{n-1},$$
	since
	$$||\hat{B}(\phi(J))||_2=c^{n-1}.$$
	
	To prove uniqueness, suppose the extremal $B(z)$ takes the form 
	$$B(z)=\prod_{j=1}^{n-1}\dfrac{z-\alpha_j}{1-\bar{\alpha_j}z},~\alpha_j\in\mathbb{C},\text{ and }|\alpha_j|\le 1,$$
	then
	$$B(J)=p_0I+p_1J+\cdots+p_{n-1}J^{n-1},$$
	or
	$$P=B(J)=\left[\begin{matrix}
	p_0 &p_1 &\cdots &p_{n-2} & p_{n-1} \\
	0 &p_0 &p_1 &\cdots &p_{n-2}\\
	0 &0 &\ddots &\ddots &\vdots\\
	0 &0 &\cdots &p_0 &p_1 \\
	0 &0 &\cdots &0 &p_0 \\
	\end{matrix}\right],$$ 
	since $J^{n}=0$. Here $p_j$ is $1 / j!$ times the $j^{th}$ derivative of $B(z)$ and can be expressed in terms of $\alpha_j$'s:
	$$p_0=(-1)^{n-1}\prod_{j=1}^{n-1}\alpha_j,$$
	$$p_1=(-1)^{n-2}\sum_{j=1}^{n-1}\left(\left(\prod_{k\ne j}\alpha_k\right)(1-|\alpha_j|^2)\right),$$
	$$\vdots$$
	We know that the $p_j$'s satisfy
	$$||P||_2\le 1,$$
	and
	\begin{equation} \label{EQUAL}
		||B(\phi(J))||_2=||D^{-1}PD||_2\le||D^{-1}||_2||P||_2||D||_2\le c^{n-1},
	\end{equation}
	with equality in both places if $B$ is extremal.
	Suppose that the first right and left singular vectors of $D^{-1}PD$ are $v_1$ and $u_1$ respectively, i.e., 
	$$u_1^*D^{-1}PDv_1=c^{n-1},~||v_1||_2=||u_1||_2=1.$$
	Then equality holds in (\ref{EQUAL}) only when $||P||_2=1$, $u_1$ is parallel to the first right singular vector of $D^{-T}=D^{-1}$, that is $e_1$ ( $e_j$'s are the standard basis), $v_1$ is parallel to the first right singular of $D$ (i.e., $e_n$), $D^{-1}u_1$ is parallel to the first left singular vector of $P$, and $Dv_1$ is parallel to the first right singular vector of $P$. Thus, we have
	$$|p_{n-1}|=1,$$
	and
	$$||[p_0~p_1\cdots p_{n-1}]^T||_2=\left(\displaystyle\sum_{j=0}^{n-1}|p_j|^2\right)^{\frac{1}{2}}\le||P||_2=1$$
	gives us
	$$p_{j}=0,~j=0,1,\cdots,n-2.$$
	Further, $p_0=0$ means $B(0)=0$, from which we get some $\alpha_{j}=0$. Without loss of generality, take $\alpha_1=0$. Then 
	$$B(z)=z\prod_{j=2}^{n-1}\dfrac{z-\alpha_j}{1-\bar{\alpha_j}z},~\alpha_j\in\mathbb{C},\text{ and }|\alpha_j|\le 1,$$
	differentiating $B(z)$ we can see $B'(0) = 0$ (i.e., $p_1=0$), which implies that one of the remaining ${\alpha}$'s is 0, i.e.,
	$$B(z)=z^2\prod_{j=3}^{n-1}\dfrac{z-\alpha_j}{1-\bar{\alpha_j}z},~\alpha_j\in\mathbb{C},\text{ and }|\alpha_j|\le 1.$$
	Do this step by step, and finally we get
	$$\alpha_1=\alpha_2=\cdots=\alpha_{n-1}=0.$$
	That is, $\hat{B}(z)=z^{n-1}$ is the only extremal Blaschke product.
	
	\subsection{3 by 3 elliptical case with one parameter}
	Since $W(J_3)$ is a disk, we can construct a matrix $A$ 
	$$A = J_{3}+bJ_{3}^*,$$
	and from the previous discussion $W(A)$ is an ellipse. C. Glader, M. Kurula, and M. Lindström argued that Courzeix's conjecture holds in this case\cite{doi:10.1137/17M1110663}. We give the exact upper bound and prove uniqueness here. 
	
	Scale $A$ by the factor $\dfrac{1}{2\sqrt{b}\cos(\frac{\pi}{4})}$, thus
	$$A=\dfrac{1}{\sqrt{2b}}\left[\begin{matrix}
	0&1&0\\ b&0&1\\ 0&b &0
	\end{matrix}\right],~ b\in(0,1].$$
	The eigenvalues of $A$ are 0 and the two foci of $W(A)$, $\pm 1$. The conformal mapping $\phi$ is the same as equation (\ref{confm}), i.e.,
	$$\phi(A)=cA,$$
	where $c=\phi(1)=k^{1/2}$. We can rewrite
	$$\phi(A)=\left[\begin{matrix}
	0 &\sqrt{\dfrac{k}{2b}} &0\\
	\sqrt{\dfrac{kb}{2}} &0 &\sqrt{\dfrac{k}{2b}}\\
	0 &\sqrt{\dfrac{kb}{2}} &0
	\end{matrix}\right]=V^{-1}TV,$$
	where $$V = \left[\begin{matrix}
	1 &0 &0\\
	0 &\dfrac{1}{\sqrt{t}} &0\\
	0 &0 &\dfrac{1}{t}
	\end{matrix}\right],$$
	$$T=\left[\begin{matrix}
	0 &\sqrt{\dfrac{kt}{2b}} &0\\
	\sqrt{\dfrac{kb}{2t}} &0 &\sqrt{\dfrac{kt}{2b}}\\
	0 &\sqrt{\dfrac{kb}{2t}} &0
	\end{matrix}\right],$$
	and
	$$t=\dfrac{b(1+\sqrt{1-k^2})}{k}.$$
	We can check that
	$$||T||_2=1,$$
	and
	$$\dfrac{1}{2}< t\le 1,\text{ if }b\in(0,1].$$
	When $b\rightarrow 0$, $t\rightarrow \dfrac{1}{2}$, i.e., $b=0$ is a removable singularity of $t$.
	
	By von Neumann's inequality, for any Blaschke product we have
	$$||B(\phi(A))||_2=||V^{-1} B(T)V||_2\le ||V^{-1}||_2||B(T)||_2||V||_2\le 1\cdot 1\cdot \dfrac{1}{t}=\dfrac{k}{b(1+\sqrt{1-k^2})}.$$
	
	Taking
	$$\hat{B}(z)=\dfrac{z^2-\alpha^2}{1-\bar{\alpha}^2z^2},$$
	where
	$$\alpha=\phi(\frac{1}{\sqrt{2}})=\sqrt{\frac{k}{1+\sqrt{1-k^2}}},$$
	we can achieve the extremal Blaschke product. That is, for this class of matrix
	$$\hat{B}(T)=\left[\begin{matrix}
	0 & 0 &1 \\ 0 &\dfrac{k}{1+\sqrt{1-k^2}} &0\\
	\dfrac{1-\sqrt{1-k^2}}{1+\sqrt{1-k^2}} &0&0
	\end{matrix}\right],$$
	$$\hat{B}(\phi(A))=\left[\begin{matrix}
	0 &0 &\dfrac{1}{t}\\ 0 &\dfrac{k}{1+\sqrt{1-k^2}} &0\\
	\dfrac{1-\sqrt{1-k^2}}{1+\sqrt{1-k^2}}t &0&0
	\end{matrix}\right],$$
	$$||\hat{B}(\phi(A))||_2=\dfrac{1}{t}=\dfrac{k}{b(1+\sqrt{1-k^2})}.$$
	
	\begin{figure}
		\centering
		\includegraphics[width=12cm]{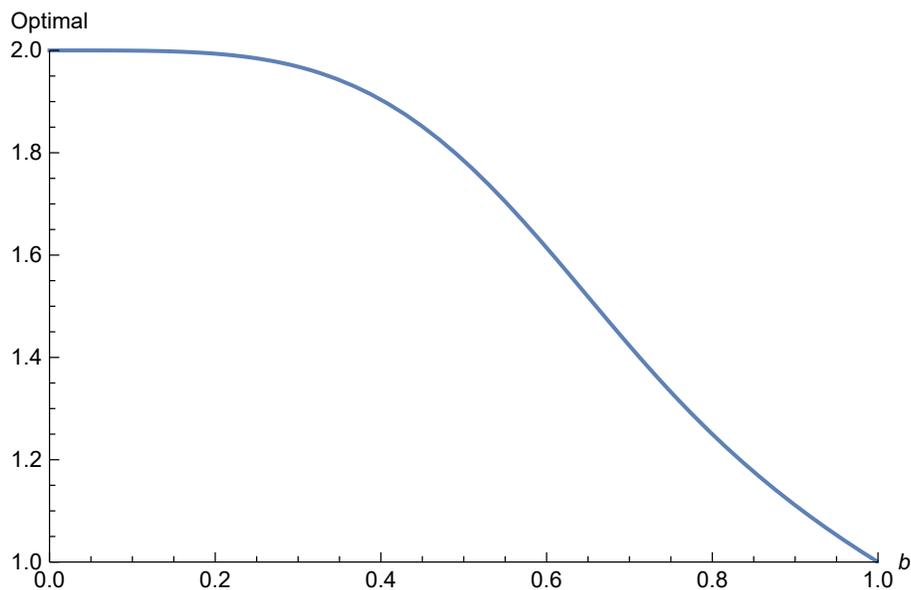}
		\caption{$||\hat{B}(\phi(A))||_2=\frac{1}{t}$ vs. $b$}
		\label{fig2.2.2}
	\end{figure}
	
	To prove uniqueness, assume the extremal Blaschke product is in the form
	$$B(z)=\dfrac{z-\alpha_1}{1-\bar{\alpha_1}z}\dfrac{z-\alpha_2}{1-\bar{\alpha_2}z}.$$
	Since the minimum polynomial of $T$ is 
	$$T^3-kT=0,$$
	we have
	$$B(T)=p_0I+p_1T+p_2T^2,$$
	where
	\begin{equation}\label{p0p1p2}
	\left\{\begin{aligned}
	p_0&=\alpha_1\alpha_2 \\
	p_1&=-\alpha_1\dfrac{1-|\alpha_2|^2}{1-\bar{\alpha_2}^2k}-\alpha_2\dfrac{1-|\alpha_1|^2}{1-\bar{\alpha_1}^2k}+k(\bar{\alpha_1}+\bar{\alpha_2})\dfrac{1-|\alpha_1|^2}{1-\bar{\alpha_1}^2k}\dfrac{1-|\alpha_2|^2}{1-\bar{\alpha_2}^2k}\\
	p_2&= -\alpha_1\bar{\alpha_2}\dfrac{1-|\alpha_2|^2}{1-\bar{\alpha_2}^2k}-\bar{\alpha_1}\alpha_2\dfrac{1-|\alpha_1|^2}{1-\bar{\alpha_1}^2k}+(1+k\bar{\alpha_1}\bar{\alpha_2})\dfrac{1-|\alpha_1|^2}{1-\bar{\alpha_1}^2k}\dfrac{1-|\alpha_2|^2}{1-\bar{\alpha_2}^2k}
	\end{aligned}\right. .
	\end{equation}
	On the other hand,
	$$||B(\phi(A))||_2\le ||V^{-1}||_2||B(T)||_2||V||_2\le\dfrac{1}{t},$$
	with equality in both places if and only if $B$ is extremal. Then, we must have
	$$||B(T)||_2=1,$$
	and $e_1$ and $e_3$ must be parallel to the first left and right singular vectors of $B(\phi(A))$ (or $B(T)$), respectively. Thus, we have
	$$B(T)=\left[\begin{matrix}
	0 &0 &x_{13} \\ x_{21} &x_{22} &0\\ x_{31} &x_{32} &0
	\end{matrix}\right],$$
	and
	$$1=|x_{13}|=||B(T)||_2.$$
	Without loss of generality, taking $x_{13}=1$, we can obtain
	$$\left\{\begin{aligned}
	p_0&=-\dfrac{k}{1+\sqrt{1-k^2}} \\
	p_1&=0 \\
	p_2&= \dfrac{2}{1+\sqrt{1-k^2}}\end{aligned}\right. ,$$
	and
	$$B(T)=\hat{B}(T).$$
	Combining with expressions (\ref{p0p1p2}) for $p_0$ and $p_1$ in terms of $\alpha_1$ and $\alpha_2$, we have
	$$\alpha_1=-\alpha_2=\sqrt{\dfrac{k}{1+\sqrt{1-k^2}}}\text{ or }\alpha_1=-\alpha_2=-\sqrt{\dfrac{k}{1+\sqrt{1-k^2}}}.$$
	One can double check that the expression (\ref{p0p1p2}) for $p_2$
	$$\dfrac{2}{1+\sqrt{1-k^2}}=p_2=-\alpha_1\bar{\alpha_2}\dfrac{1-|\alpha_2|^2}{1-\bar{\alpha_2}^2k}-\bar{\alpha_1}\alpha_2\dfrac{1-|\alpha_1|^2}{1-\bar{\alpha_1}^2k}+(1+k\bar{\alpha_1}\bar{\alpha_2})\dfrac{1-|\alpha_1|^2}{1-\bar{\alpha_1}^2k}\dfrac{1-|\alpha_2|^2}{1-\bar{\alpha_2}^2k}$$
	is also satisfied. In sum, the extremal Blaschke product is unique
	$$B(z)=\hat{B}(z).$$
	
	Further, we can find the relation of Figure \ref{fig2.2.2} with Figure \ref{fig1.3.1} by showing the equality of $k(q^2)$ and $k(q)$
	$$\sqrt{k(q^2)}=\dfrac{k(q)}{1+\sqrt{1-k(q)^2}},~q\in [0,1].$$
	\begin{proof}
		We know that the Jacobi theta functions satisfy
		\begin{equation} \label{JacobiTheta}
		\theta_3(q)^4 = \theta_2(q)^4+\theta_4(q)^4,
		\end{equation}
		and we can rewrite them in terms of the Dedekind eta function
		$$\theta_2(q)=\dfrac{2\eta(2\tau)^2}{\eta(\tau)},$$
		$$\theta_3(q)=\dfrac{\eta(\tau)^5}{\eta(2\tau)^2\eta(\frac{\tau}{2})^2},$$
		$$\theta_4(q)=\dfrac{\eta(\frac{\tau}{2})^2}{\eta(\tau)},$$
		where $q=\exp(\pi i \tau)$.
		Substituting all $\theta_j$ into equality (\ref{JacobiTheta}), we find
		$$\eta(\tau)^{24}=\eta(\frac{\tau}{2})^8\eta(2\tau)^8(16\eta(2\tau)^8+\eta(\frac{\tau}{2})^8).$$
		It follows that
		$$\eta(\frac{\tau}{2})^8=\dfrac{-8\eta(2\tau)^{12}+\sqrt{64\eta(2\tau)^{24}+\eta(\tau)^{24} }}{\eta(2\tau)^4},$$
		$$\eta(2\tau)^8=\dfrac{-\eta(\frac{\tau}{2})^{12}+\sqrt{64\eta(\tau)^{24}+\eta(\frac{\tau}{2})^{24} }}{32\eta(\frac{\tau}{2})^4},$$
		or
		$$\eta(4\tau)^8=\dfrac{-\eta(\tau)^{12}+\sqrt{64\eta(2\tau)^{24}+\eta(\tau)^{24} }}{32\eta(\tau)^4}.$$
		Expressing $\eta(4\tau)$ and $\eta(\frac{\tau}{2})$ in terms of $\eta(\tau)$ and $\eta(2\tau)$ on the left hand side, we have an identity
		$$\eta(\tau)^2\eta(4\tau)^4\left(\eta(\tau)^{12}+\eta(2\tau)^4\eta(\frac{\tau}{2})^8\right)=2\eta(2\tau)^{14}\eta(\frac{\tau}{2})^4.$$
		This can be written in terms of $\theta_j$,
		$$\theta_2(q^2)\left(\theta_3(q)^2+\theta_4(q)^2\right)=\theta_2(q)^2\theta_3(q^2),$$
		which is equivalent to the equality 
		$$\sqrt{k(q^2)}=\dfrac{k(q)}{1+\sqrt{1-k(q)^2}}.$$
	\end{proof}
	
	Figure \ref{fig2.2.2} and Figure \ref{fig1.3.1} match the functions $\sqrt{\dfrac{k(b^4)}{b^2}}$ and $\sqrt{\dfrac{k(b^2)}{b}}$ respectively.
	
	\section{One Example of Nonuniqueness}
	Let
	$$A=\left[\begin{matrix}
	0 &1 &0\\ 0 &0 &1-t\\0 &0 &0
	\end{matrix}\right],~t\in[0,\sqrt{3}-1].$$
	$W(A)$ is a disk centered at the origin with radius 
	$$r=\dfrac{\sqrt{(1+(1-t)^2)}}{2}\le 1.$$
	The conformal mapping $\phi$ that maps $W(A)$ to the unit disk $\mathbb{D}$ is
	$$\phi(z)=\dfrac{1}{r}z,$$
	which gives 
	$$\phi(A)=\dfrac{1}{r}A=\left[\begin{matrix}
	0 &\frac{1}{r} &0\\0 &0 &\frac{1-t}{r} \\ 0&0 &0
	\end{matrix}\right]=D\left[\begin{matrix}
	0 &1 &0\\ 0&0 &1 \\0 &0 &0
	\end{matrix}\right]D^{-1}=DJ_3D^{-1},$$
	where
	$$D=\left[\begin{matrix}
	1 &0 &0\\ 0&r &0\\ 0&0&\dfrac{r^2}{1-t}
	\end{matrix}\right], ~J_3=\left[\begin{matrix}
	0 &1 &0\\ 0&0 &1 \\0 &0 &0
	\end{matrix}\right].$$
	
	Since $0\le t\le\sqrt{3}-1$, we get
	$$\dfrac{\sqrt{3}-1}{\sqrt{2}}\le r\le \dfrac{1}{\sqrt{2}},$$
	and
	$$\dfrac{1}{2}\le \dfrac{r^2}{1-t}\le 1.$$
	Then, for any Blaschke product $B$ we can obtain
	$$B(\phi(A))=DB(J_3)D^{-1},$$
	and
	$$||B(\phi(A))||_2\le||D||_2||B(J_3)||_2||D^{-1}||_2=1\times 1\times ||D^{-1}||_2=\max(\dfrac{1-t}{r^2},\dfrac{1}{r}).$$
	Taking $B_1(z)=z$, we find
	$$||B_1(\phi(A))||_2=\dfrac{1}{r},$$
	and taking $B_2(z)=z^2$, we find
	$$||B_2(\phi(A))||_2=\dfrac{1-t}{r^2}.$$
	In addition, on the interval $[0,\sqrt{3}-1]$, the function
	$h(t)=\dfrac{1-t}{r}=\dfrac{2}{\sqrt{1+\frac{1}{(1-t)^2}}}$
	decreases monotonically as $t$ increases, and we find the point $t_0$ where $||B_1(\phi(A))||_2$ starts to exceed $||B_2(\phi(A))||_2$ by solving
	$$\dfrac{1}{r}=\dfrac{1-t_0}{r^2}$$
	or
	$$h(t_0)=1,$$
	which gives us
	$$t_0=1-\dfrac{1}{\sqrt{3}}.$$
	When $t\le t_0$, the extremal Blaschke product is $z^2$; when $t\ge t_0$, it is $z$; at the point $t_0$, both Blaschke products ($z$ and $z^2$) give us the same extremal result
	$$||B(\phi(A))||_2=\sqrt{3}.$$
	\begin{figure}[htbp]
		\centering
		\includegraphics[width=12cm]{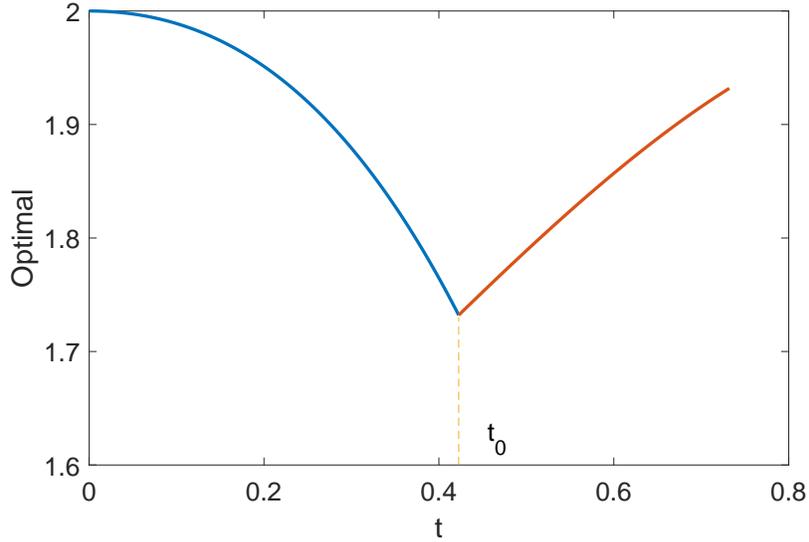}
		\caption{$||\hat{B}(\phi(J_t))||_2$ vs. $t$. Blue is $\hat{B}(z)=z^2$, and red is $\hat{B}(z)=z$.}
		\label{fig2.2.1}
	\end{figure}
	Figure \ref{fig2.2.1} shows a plot of $\max_{B}||B(\phi(A))||_2=\max{||B_1{\phi(A)}||_2,||B_2{\phi(A)}||_2}$ as $t$ ranges between $0$ and $\sqrt{3}-1$. In addition to showing nonuniqueness, the extremal Blaschke product has degree less than the maximal degree $n-1$. This has been observed in numerical experiments, and we have now proved it for this class of matrices when $1-1/\sqrt{3} <t<\sqrt{3}-1$.
	
	At the point $t=t_0=1-\frac{1}{\sqrt{3}}$, we have 
	$$\phi(A) = \left[\begin{matrix}
	0 & \sqrt{3} & 0\\
	0 & 0 & 1\\
	0 & 0 & 0
	\end{matrix}\right]\text{ and }
	D=\left[\begin{matrix}
	1 & 0 & 0\\
	0 & \frac{1}{\sqrt{3}} & 0\\
	0 & 0 & \frac{1}{\sqrt{3}}
	\end{matrix}\right].$$
	Suppose that the first right and left singular vectors of $DB(J_3)D^{-1}$ are $v_1$ and $u_1$ respectively, i.e., 
	$$u_1^*DB(J_3)D^{-1}v_1=\sqrt{3},~||v_1||_2=||u_1||_2=1.$$
	We can see that $u_1$ must be parallel to the first right singular vector of $D^{T}=D$ (i.e., $e_1$), $v_1$ should be parallel to the first right singular vector of $D^{-1}$ (i.e., $e_2$ or $e_3$). In fact, the choices of $v_1$ ($v_1=e_2$ or $v_1=e_3$) correspond to different extremal Blachke products $z$ and $z^2$, respectively. 
	
	\section{Conclusion}
	In the discussion above, we focused on decomposing the conformal mapping of matrix $A$, 
	$$\phi (A) = D C D^{-1},$$ 
	where the condition number of $D$, $\kappa(D)$, is no more than 2 and $C$ is a contraction ($||C||_2\le 1$), and using von Neumann's inequality to argue that for any Blaschke product $B$, we have 
	$$\| B( \phi (A)) \|_2 \leq
	\kappa (D)\le 2,$$ 
	since $\| B(C) \|_2 \leq 1$. The construction of matrix $D$ is the key to argue the exact bound, the extremal Blaschke product and its uniqueness. We may extend 
	this decomposition to more general cases in the future.
	
	\bibliographystyle{plain}
	\bibliography{bibtxt} 

\begin{thebibliography}{1}

\bibitem{Discussions}
Discussions at {AIM} workshop of {Crouzeix's conjecture}.

\bibitem{caldwell2017extensions}
Trevor Caldwell, Anne Greenbaum, and Kenan Li.
\newblock Some extensions of the crouzeix--palencia result.
\newblock {\em SIAM J. Matrix Anal. Appl.}, 39(2):769--780, May 2018.

\bibitem{doi:10.1137/17M1116672}
M.~Crouzeix and C.~Palencia.
\newblock The numerical range is a $(1+\sqrt{2})$-spectral set.
\newblock {\em SIAM Journal on Matrix Analysis and Applications},
  38(2):649--655, 2017.

\bibitem{Crouzeix2004}
Michel Crouzeix.
\newblock Bounds for analytical functions of matrices.
\newblock {\em Integral Equations and Operator Theory}, 48(4):461--477, Apr
  2004.

\bibitem{CROUZEIX2007668}
Michel Crouzeix.
\newblock Numerical range and functional calculus in hilbert space.
\newblock {\em Journal of Functional Analysis}, 244(2):668 -- 690, 2007.

\bibitem{doi:10.1137/17M1110663}
C.~Glader, M.~Kurula, and M.~Lindström.
\newblock Crouzeix's conjecture holds for tridiagonal 3 x 3 matrices with
  elliptic numerical range centered at an eigenvalue.
\newblock {\em SIAM Journal on Matrix Analysis and Applications},
  39(1):346--364, 2018.

\bibitem{GREENBAUM2018225}
Anne Greenbaum and Michael~L. Overton.
\newblock Numerical investigation of crouzeix's conjecture.
\newblock {\em Linear Algebra and its Applications}, 542:225 -- 245, 2018.
\newblock Proceedings of the 20th ILAS Conference, Leuven, Belgium 2016.

\bibitem{Neumann}
Johann~Von Neumann.
\newblock Eine spektraltheorie für allgemeine operatoren eines unitären
  raumes. erhard schmidt zum 75. geburtstag in verehrung gewidmet.
\newblock {\em Mathematische Nachrichten}, 4(1-6):258--281, 1950.

\bibitem{Szego}
Gabor Szego.
\newblock Conformal mapping of the interior of an ellipse onto a circle.
\newblock {\em The American Mathematical Monthly}, 57(7):474--478, 1950.

\end{thebibliography}
\end{document}